\documentstyle{amsppt}
\magnification=\magstep1
\topmatter
\title 
On Vorontsov's theorem on K3 surfaces
\endtitle
\author Keiji Oguiso and De-Qi Zhang
\endauthor
\abstract
Let  $X$  be a K3 surface with the N\'eron-Severi lattice
$S_X$  and transcendental lattice  $T_X$.
Nukulin considered the kernel  $H_X$ of the natural representation
$\text{Aut}(X) \longrightarrow O(S_X)$  and proved that $H_{X}$ is a finite 
cyclic group with $\varphi(h(X))) | t(X)$  and acts 
faithfully on the space $H^{2,0}(X) = {\bold C}\omega_{X}$, where 
$h(X) =$ ord$(H_X)$, $t(X) =$ rank $T_X$  and
$\varphi$  is the Euler function.
Consider the extremal case where  $\varphi(h(X)) = t(X)$.
In the situation where  $T_{X}$ is unimodular,
Kondo has determined the list of  $t(X)$, as well as the actual realizations,
and showed that  $t(X)$ alone uniquely determines the isomorphism class
of  $X$  (with  $\varphi(h(X)) = t(X)$).
We settle the remaining situation where  $T_X$  is not unimodular.  
Together, we provide the proof for the theorem announced by Vorontsov.
\endabstract
\endtopmatter

\document
{\bf Introduction}

\par \vskip 0.5pc
A K3 surface is, by definition, a simply connected
smooth projective surface over the complex numbers  ${\bold C}$  
with a nowhere vanishing holomorphic $2-$form.
For a K3 surface $X$, we denote by
$S_X$, $T_X$  and  $\omega_X$  the N\'eron-Severi lattice,
the transcendental lattice and 
a nowhere vanishing holomorphic $2-$form of  $X$.
We write $t(X) = \text{rank} T_X$.

\par
Nukulin [Ni1] considered the kernel  $H_X$ of the natural representation
$\text{Aut}(X) \longrightarrow O(S_X)$  and proved that $H_{X}$ is a finite 
cyclic group with $\varphi(\text{ord}(H_{X})) | t(X)$ and acts faithfully on the space $H^{2,0}(X) = {\bold C}\omega_{X}$, where $\varphi$  is the Euler function. We set $h(X) = \text{ord}(H_{X})$. 
The interesting case here is when  $\varphi(h(X)) = t(X)$.

\par
Kondo [Ko, main Theorem] has studied the case where $T_{X}$ is unimodular and 
shown the following complete classification:

\par \vskip 1pc
{\bf Theorem 1.} {\it Set  $\Sigma := \{66, 44, 42, 36, 28, 12 \}$.}

\par
(1) {\it Let $X$ be a  $K3$  surface with $\varphi(h(X)) = t(X)$ whose 
transcendental lattice $T_{X}$ is unimodular. Then $h(X) \in \Sigma$.}

\par
(2) {\it Conversely, for each  $N \in \Sigma$,
there exists, modulo isomorphisms, a unique
$K3$  surface  $X$  such that  $h(X) = N, \varphi(h(X)) = t(X)$. Moreover, 
$T_X$ is unimodular for this $X$.}

\par \vskip 1pc
In the case where $T_{X}$ is not unimodular, about 15 years ago, 
Vorontsov [Vo] announced the following complete classification:

\par \vskip 1pc
{\bf Theorem 2.} {\it Set  $\Omega := \{3^{k} (1 \leq k \leq 3), 
5^{l} (l = 1,2), 7, 11, 13, 17, 19\}$.}

\par 
(1) {\it Assume that  $X$  is a  $K3$  surface
satisfying  $\varphi(h(X)) = t(X)$ 
and that  $T_X$  is non-unimodular.  Then  $h(X) \in \Omega$.}

\par
(2) {\it Conversely, for each  $N \in \Omega$,
there exists, modulo isomorphisms, a unique
$K3$  surface  $X$  such that  $h(X) = N, \varphi(h(X)) = t(X)$. Moreover, 
$T_X$  is non-unimodular for this $X$.}

\par \vskip 1pc
However, till now, he gave neither proof of this Theorem
nor construction of such K3 surfaces. In fact the original statement 
of (1) in [V] was weaker than here.

\par 
Later, in [Ko, Sections 6 and 7], Kondo has sharpened the statement of (1) 
as in the present form 
and also given a complete proof of the statement (1). 
He has also shown the existence part of the statement (2) 
by constructing such K3 surfaces explicitly as follows: 

\par \vskip 1pc
{\bf Kondo's Example}

\par \vskip 0.5pc
For each $h \in \Omega$, the following pair $(X_{h}, <g_{h}>)$ 
of a K3 surface $X_{h}$ defined by the indicated Weierstrass equation 
(or a weighted homogeneous equation in a weighted projective space)
and its cyclic automorphism group $<g_{h}>$ satisfies 
$H_{X_{h}} = <g_{h}>$ and  $\varphi(h(X_{h})) = t(X_{h})$ and $T_{X}$ 
is not unimodular: 

\par \vskip 0.5pc
$X_{19} : y^2 = x^3 + t^{7}x + t, \,\,\,$
$g_{19}^*(x,y,t) = (\zeta_{19}^{7}x, \zeta_{19}y, \zeta_{19}^{2}t)$;

\par
$X_{17} : y^2 = x^3 + t^{7}x + t^{2}, \,\,\,$
$g_{17}^*(x,y,t) = (\zeta_{17}^{7}x, \zeta_{17}^{2}y, \zeta_{17}^{2}t)$;

\par
$X_{13} : y^2 = x^3 + t^{5}x + t^{4}, \,\,\,$ 
$g_{13}^*(x,y,t) = (\zeta_{13}^{5}x, \zeta_{13}y, \zeta_{13}^{2}t)$; 

\par
$X_{11} : y^2 = x^3 + t^{5}x + t^{2}, \,\,\,$ 
$g_{11}^*(x,y,t) = (\zeta_{11}^{5}x, \zeta_{11}^{2}y, \zeta_{11}^{2}t)$; 

\par
$X_{7} : y^2 = x^3 + t^{3}x + t^{8}, \,\,\,$ 
$g_{7}^*(x,y,t) = (\zeta_{7}^{3}x, \zeta_{7}y, \zeta_{7}^{2}t)$;  

\par
$X_{25} : \{y^{2} + x_{0}^{6} + x_{0}x_{1}^{5} + x_{1}x_{2}^{5} = 0\} \subset 
{\bold P}(1,1,1,3)$,
\par \hskip 2pc
$g_{25}^{*}([x_{0} : x_{1} : x_{2} : y])
= [x_{0} : \zeta_{25}^{20}x_{1} : \zeta_{25}x_{2} : y]$;

\par
$X_{5} : y^2 = x^3 + t^{3}x + t^{7}, \,\,\,$ 
$g_{5}^*(x,y,t) = (\zeta_{5}^{3}x, \zeta_{5}^{2}y, \zeta_{5}^{2}t)$;  

\par
$X_{27} : y^2 = x^3 + t(t^{9}-1), \,\,\,$  
$g_{27}^*(x,y,t) = (\zeta_{27}^{2}x, \zeta_{27}^{3}y, \zeta_{27}^{6}t)$;

\par
$X_{9} : y^2 = x^3 + t^{5}(t^{3}-1), \,\,\,$ 
$g_{9}^*(x,y,t) = (\zeta_{9}^{2}x, \zeta_{9}^{3}y, \zeta_{9}^{3}t)$; 

\par
$X_{3} : y^2 = x^3 + t^{2}(t^{10}-1), \,\,\,$ 
$g_{3}^*(x,y,t) = (\zeta_{3}x, y, t)$.

\par \vskip 1pc
However, Kondo did not touch the uniqueness part of (2), either. Only 
the uniqueness in the case where $h(X) = 5^{2}$ has been just settled 
by [MO, Theorem 3].

\par \vskip 1pc
The main purpose of this short article is to give a complete proof 
for the uniqueness part of (2) to guarantee Vorontsov's Theorem. 
This together with Kondo's Theorem completes the classification 
of K3 surfaces $X$ with $\varphi(h(X)) = t(X)$. 

\par \vskip 1pc
We shall also show the following strong uniqueness result as 
an application of Theorem 2:

\par \vskip 1pc
{\bf Corollary 3.} 
{\it Let $X$ be a  $K3$  surface with an automorphism $g$ of order 
$I \in \{19, 17, 13\}$, 
the three largest possible prime orders. Then we have:}

\par
{\it $(X, <g>) \simeq (X_{19}, <g_{19}>)$  when  $I = 19$; 
$(X, <g>) \simeq (X_{17}, <g_{17}>)$  when  $I = 17$; and
$(X, <g>) \simeq (X_{13}, <g_{13}>)$  when  $I = 13$,

\par
where $(X_{I}, <g_{I}>)$ are pairs defined in Kondo's example.}

\par \vskip 1pc
Besides its own interest, our motivation for this project lies 
also in its applicability to the study of log Enriques surfaces initiated 
by the second author ([Z1]). We should also mention here that 
log Enriques surfaces are regarded as a log version of K3 surfaces 
and play an increasingly important role in higher dimensional 
algebraic geometry.  For instance, base spaces of elliptically 
fibered Calabi-Yau threefolds $\Phi_{D} : X \rightarrow S$ 
with $D.c_{2}(X) = 0$ are necessarily log Enriques surfaces ([Og]). 

\par 
A log Enriques surface  $Z$  is, by definition, 
a projective rational surface with 
at worst quotient singularities, or in other words, 
at worst klt singularities and with numerically trivial canonical 
Weil divisor. Passing to the maximal crepant partial resolution, 
we may also assume in the definition the following
maximality for  $Z$:

\par \vskip 1pc
$(*)$ {\it any birational morphism  $Z' \rightarrow Z$
from another log Enriques surface $Z'$  must be an isomorphism.} 

\par \vskip 1pc
For a log Enriques surface $Z$, we define the canonical 
index $I(Z)$ or index for short, by 
$$
I(Z) := 
\text{\rm min} \{n \in {\bold Z}_{>0} \vert \Cal O_{Z}(nK_{Z}) 
\simeq \Cal O_{Z}\}. 
$$ 
A log Enriques surface of index  $I$  is closely 
related to a K3 surface admitting a non-symplectic group 
($\simeq {\bold Z}/I{\bold Z}$) action 
via the canonical cover and its minimal resolution:
$$
X \overset{\nu} \to{\longrightarrow} \overline{X} 
:= Spec(\oplus_{n = 0}^{I-1}\Cal O_{Z}(-iK_{Z})) 
\overset{\pi} \to{\rightarrow} Z. 
$$ 
In fact, it is well known that $\overline{X}$ is either 
an abelain surface or a normal K3 surface with at worst 
Du Val singularities and that $\pi : \overline{X} \rightarrow Z$ is 
a cyclic Galois cover of order $I$ which acts faithfully 
on the space $H^{0}({\overline{X}}, \Cal O_{\overline{X}}(K_{\overline{X}})) 
= {\bold C} \omega_{\overline{X}}$ and is ramified only over 
$\text{\rm Sing}(Z)$ ([Ka], [Z1]). 
\par 

In the case where $\overline{X}$ is an abelain surface, Blache [Bl] 
shows that there are exactly two such log Enriques surfaces 
up to isomorphisms. 
\par 
Let us consider the case where $X$ is a K3 surface. 
In [OZ1, 2, 3], we regard the rank of the sublattice of $S_{X}$ 
generated by the exceptional curves of $\pi$ as an invariant 
to measure how bad Sing$(Z)$ and classified the worst case, namely, 
'the extremal' case where the rank is 19. 
As a result, we found that there exist exactly 7 such surfaces up to 
isomorphisms. However one of them is of index 2 and the others 
are all of index 3.  Note that these indices are rather small.
\par 
Now, as a counter part, it is also interesting to consider the canonical index 
$I(Z)$ as an invariant measuring how bad  Sing$(Z)$ is. 
It is known that $2 \leq I(Z) \leq 21$ 
and $I(Z) \in \{2, 3, 5, 7, 9, 11, 13, 17, 19\}$ 
if $I(Z)$ is prime ([Z1], [Bl]). 
\par \vskip 1pc
As an application of Corollary 3, we show the following uniqueness 
result for log Enriques surfaces  $Z$  with the three largest prime indices: 

\par \vskip 1pc
{\bf Corollary 4.} 
{\it Let $Z$ be a log Enriques surface with $I(Z) = 19, 17$  or  $13$
satisfying the maximality $(*)$.
Then we have:

\par
$Z \simeq Z_{19} := X_{19}/<g_{19}>$  when  $I(Z) = 19$;

\par
$Z \simeq Z_{17} := X_{17}/<g_{17}>$  when  $I(Z) = 17$; and

\par
$Z \simeq Z_{13} := \overline{X_{13}}/<g_{13}>$  when  $I(Z) = 13$,}

\par
{\it where $(X_{I}, g_{I})$ are pairs defined in Kondo's example and 
$\overline{X_{13}}$ is the surface obtained from  $X_{13}$
by contracting the unique rational curve in the fixed 
locus  $X_{13}^{g_{13}}$.}

\par \vskip 1pc
The second author constructed log Enriques surfaces of indices $19, 17, 13$ 
in a completely different way ([Z1]). However, it looks very hard to 
show directly that they are isomorphic to $Z_{19}$, $Z_{17}$ and $Z_{13}$.   

\par \vskip 1.5pc 
{\bf Acknowledgement.} The present version of this article 
has been completed during 
the first author's stay in Singapore in March 1997 under financial 
support from the JSPS and the National University of Singapore. 
He would like to express the gratitude to both of them.  

\par \vskip 2pc
{\bf \S 1. Existence of Jacobian fiber space structures}  

\par \vskip 0.5pc
Throughout this section we assume that $X$ is a K3 surface with
$\varphi(h(X)) = t(X)$ and with $p^{r} = N = h(X) \in \Omega$,
where  $p$  is prime and fix a generator $g$ of $H_{X}$ 
with $g^{*}\omega_{X} = \zeta_{N}\omega_{X}$.
In what follows, set $S_{X}^{*} = \text{Hom}(S_{X}, {\bold Z})$, $T_{X}^{*} 
= \text{Hom}(T_{X}, {\bold Z})$ and regard $S_{X} \subset S^{*} 
\subset S_{X} \otimes {\bold Q}$, $T_{X} \subset T^{*} \subset T_{X} 
\otimes {\bold Q}$ via the bilinear form of $S_{X}$ and $T_{X}$ induced by 
the cup product on $H^{2}(X, {\bold Z})$. We denote by $l(S_{X})$ the minimal 
number of generators of the finite abelian group $S_{X}^{*}/S_{X}$. 
We call $S_{X}$ $p-$elementary if there exists a non-negative integer $a$ 
such that $S_{X}^{*}/S_{X}$ is isomorphic to $({\bold Z}/p)^{\oplus a}$. 
In this case we denote this $a$ by $l(S_{X})$.  
\par 
Recall that $S_{X}$ (resp. $T_{X}$) is an even lattice of signature 
$(1, \text{rank}S_{X} - 1)$ (resp. of signature 
$(2, \text{rank}T_{X} - 2)$) and $\text{rank}S_{X} + \text{rank}T_{X} = 22$. 

\par \vskip 1pc
The goal of this section is to show the following: 

\par \vskip 1pc
{\bf Proposition(1.1).} 
{\it $X$ admits a Jacobian fibration $\Phi : X \rightarrow {\bold P}^{1}$.} 
if $N \not= 25$.

\par \vskip 1pc
First we notice the following: 

\par \vskip 1pc
{\bf Lemma(1.2)} ([MO], [Ni 1]).

\par 
(1) {\it Each eigenvalue of $g^{*}\vert T_{X}$ is a primitive 
$N-$th root of  $1$.}

\par
(2) {\it $\text{Ann}(T_{X}) = <\Phi_{N}(g^*)>$, and 
$T_{X}$ is then naturally a torsion free 
${\bold Z}[<g^*>] /<\Phi_{N}(g^*)>-$module, 
where  $\Phi_N(x)$  denotes the minimal polynomial
over  ${\bold Q}$  of a primitive  $N$-th root
of  $1$.}

\par
(3) {\it under the identification 
${\bold Z}[<g^*>] /<\Phi_{N}(g^*)> = 
{\bold Z}[\zeta_{N}]$ through the correspondence 
$g^* (\text{mod} <\Phi_{N}(g^*)>) \leftrightarrow \zeta_{N}$, 
$T_{X} \simeq {\bold Z}[\zeta_{N}]$ as 
${\bold Z}[\zeta_{N}]-$modules.}

\par \vskip 1pc
{\it Proof.} This is proved in [MO, Lemma(1.1)]. But it is so easy that we reproduce 
the verfication here from [MO]. 
The statement (1) is shown by Nukulin ([Ni1, Theorem 3.1, Cor. 3.3]). 
The satement (2) is 
a simple reinterpretation of (1) in terms of group algebra. 
Recall that torsion free modules are in fact free if the coefficient ring is 
PID. Now, combining (2) with the fact that ${\bold Z}[\zeta_{N}]$ 
is PID for $N \in \Omega$ [MM, main Theorem], we get the assertion (3). 

\par \vskip 1pc
{\bf Lemma(1.3).} {\it $S_{X}$ is a $p-$elementary 
lattice with  $l(S_{X}) = 1$.}

\par \vskip 1.5pc
{\it Proof.} Since there exists a natural isomorphism 
$T_{X}^{*}/T_{X} \simeq  S_{X}^{*}/S_{X}$ which commutes 
with the action of $\text{Aut}(X)$, it is enough to show that 
$T_{X}^{*}/T_{X} \simeq {\bold Z}/p$.
Since $g^{*} \vert S_{X} = id$ by the definition of $H_{X}$, 
$g^{*} \vert (S_{X}^{*}/S_{X}) = id$, whence 
$g^{*} \vert (T_{X}^{*}/T_{X}) = id$. 
This means $g^{*}(x) \equiv x$ (mod $T_X$) 
for each $x \in T_{X}^{*}$. Set $n = N/p$ and 
$h = g^{n}$. Then $h$ is of order $p$. Using (1.2)(1), we get 
$px \equiv x + h^{*}(x) + ... + (h^{*})^{p-1}(x) = (1 + h^{*} + ... 
+ (h^{*})^{p-1})(x) = 0 (\text{mod} T_{X})$. Thus, $T_{X}^{*}/T_{X}$ 
is $p-$elementary. We determine $l(T_{X})$. 

\par 
We shall treat the case where $N = p$. 
The verification for the case where  $N = 3^{2}, 3^{3}, 5^{2}$ 
is quite similar and left to the reader as an exercise. 
(cf. [MO, Claim(3.4)] for the case where $N = 5^{2}$). 
Let $e_{i}$ ($i = 1, ..., p-1$) be a ${\bold Z}-$basis of $T_{X}$ 
corresponding to the ${\bold Z}-$basis 
$1, \zeta_{p}, ..., \zeta_{p}^{p-2}$ 
of ${\bold Z}[\zeta_{p}]$ via the isomorphism in (1.2). 
Then $g^{*}(e_{i}) = e_{i+1}$ for $i = 1 ,..., p-1$ and 
$g^{*}(e_{p-1}) = -(e_{1} + e_{2} + ... + e_{p-1})$ 
(corresponding to the equality 
$\Phi_{p}(\zeta_{p}) = 0$ in ${\bold Z}[\zeta_{p}]$).

\par
Choose $y \in T_{X}^{*} (\subset T_{X}\otimes {\bold Q})$ arbitrary. 
Since $T_{X}^{*}/T_{X}$ is $p-$elementary, 
we can write $y = 1/p (\sum_{i=1}^{p-1} a_{i}e_{i})$, 
where $a_{i} \in {\bold Z}$.
Then  $g^{*}(y) - y = 1/p(-(a_{1}+a_{p-1})e_{1} + 
\sum_{i=1}^{p-3}(a_{i} -a_{p-1} -a_{i+1})e_{i+1} + 
(a_{p-2}-2a_{p-1})e_{p-1})$. 
Since $g^{*}\vert (T_{X}^{*}/T_{X}) = id$,
we have $g^{*}(y) -y \in T_{X}$,
whence $a_{1}+a_{p-1} \equiv 0$, $a_{i} -a_{p-1} -a_{i+1} \equiv 0$ and 
$a_{p-2}-2a_{p-1} \equiv 0$ (mod $p$). 
This implies $a_{i} \equiv ia_{1}$ (mod $p$) and then
$y = a_{1}\times (1/p)(e_{1} + 2e_{2} + ... + (p-1)e_{p-1})$ 
in $T_{X}^{*}/ T_{X}$. Thus, 
$T_{X}^{*}/ T_{X} = < (1/p)(e_{1} + 2e_{2} + ... + (p-1)e_{p-1})> 
\simeq {\bold Z}/p$  because  $l(T_{X}) \not= 0$ if $N \in \Omega$
(cf. [Ko]).  This implies the result. 

\par \vskip 1pc
{\it Proof of Proposition }(1.1).
Let $U$ be the even unimodular hyperbolic lattice of rank $2$. 
If $N \in \Omega - \{5^{2}\}$, then $\text{rank}(S_{X}) 
\geq 4 =  3 + l(S_{X})$. 
We can then apply the so-called splitting Theorem due to 
Nikulin [Ni3, Cor. 1.13.5] 
for $S_{X}$ to split $U$ out from $S_{X}$, namely, $S_{X} \simeq U \oplus S'$. 
Now the result follows from [Ko, Lemma 2.1]. Q.E.D.

\par \vskip 2pc
{\bf \S 2. Uniqueness theorem when  $h(X) = 3^{3}, 3^{2}, 3$} 

\par \vskip 0.5pc
In this section we show the uniqueness of K3 surfaces $X$ with 
$\varphi(h(X)) = t(X)$ and with 
$N := h(X) = 3^{3} (\text{resp.} 3^{2},\text{resp.} 3)$. 
Let us set $H_{X} = <g>$.
Then $\text{rank}S_{X} = 22 - t(X) = 4$, (resp. $16$, resp. $20$). 
Since $S_{X}$ is an even hyperbolic $3-$elementary lattice with 
$\ell(S_{X}) = 1$ by (1.3), 
applying [RS, Section 1], we find that  
$S_{X} \simeq U \oplus A_{2}$, $U \oplus E_{8} \oplus E_{6}$, and 
$U \oplus E_{8} \oplus E_{8} \oplus A_{2}$. 
Thus $X$ has a Jacobian fibration 
$\Phi : X \rightarrow {\bold P}^{1}$ whose reducible fibers are 
exactly $I_{3}$ or $IV$, (resp. $II^{*} + IV^{*}$, resp. $II^{*} + II^{*} + I_{3}$ or $II^{*} + II^{*} + IV$). 

\par
Since $g^{*}\vert S_{X} = id$, there exists $\overline{g} \in 
\text{Aut}({\bold P}^{1})$ such that $\Phi \circ g = \overline{g} \circ \Phi$.
Note also that each smooth rational curve on $X$ must be 
$g-$stable whence each reducible fiber of $\Phi$ is also $g-$stable. 

\par
First consider the case where $N = 3$.
Since there exist three reducible fiber, 
$\overline{g} = id$. Thus each smooth fiber $E$ is $g-$stable and 
$(g\vert E)^{*}\omega_{E} = \zeta_{3}\omega_{E}$. Thus, the $J-$invariant map
$J : {\bold P}^{1} \rightarrow {\bold P}^{1}$ is constant 
$j({\bold C}/ {\bold Z} + {\bold Z} \zeta_{3}) =0$. 
In particular, each singular 
fiber is either of Type $II$, $II^{*}$, $IV$ or $IV^{*}$ by the 
classification of singular fibers ([Kd]). Thus, the reducible fibers of $\Phi$ are 
$II^{*} + II^{*} + IV$. 
We may adjust an inhomogeneous coordinate $t$ of the base so that
$X_{-1}$  and  $X_{1}$ are of type $II^{*}$ and $X_{0}$ is of type $IV$. 
Since $\chi_{top}(X) = 24 = \chi_{top}(X_{1}) + \chi_{top}(X_{-1}) + 
\chi_{top}(X_{0})$, there are no other singular fibers.

\par 
Let us determine the minimal Weierstrass equation 
$y^{2} = x^{3} + a(t)x + b(t)$ of $\Phi$. We use the notation in 
[Ne, Table in the last page].
Since $J(t) = 4a(t)^{3}/ (4a(t)^{3} + 27b(t)^{2}) = 0$, we have $a(t) = 0$ 
as polynomials. Thus, $\Delta(t) = 27b(t)^{2}$. This has exactly two 
zero's of order $10 (\text{mod} 12)$ at $t = 1, -1$ and one zero of order 
$4 (\text{mod} 12)$ at $t = 0$. Note that $\text{deg}\Delta(t) \leq 24$, 
because $X$ is a K3 surface. Thus, $\Delta(t) = C(t^{10} -1)^{2}t^{4}$ for 
some constant $C\not= 0$, whence $b(t) = c(t^{10} -1)t^{2}$ for some 
constant $c \not= 0$. This means the equation is written as 
$y^{2} = x^{3} + c(t^{10} -1)t^{2}$. 
Then changing the coordinates $x, y$  to  $c^{1/3}x, c^{1/2}y$, 
we normalise this equation as $y^{2} = x^{3} + (t^{10} -1)t^{2}$. 
This shows that $X$ is isomorphic to 
the Jacobian K3 surface $y^{2} = x^{3} + (t^{10} -1)t^{2}$. 

\par
Next consider the case where $N = 9$. 
We may take an imhomogeneous coordinate $t$  
so that  $X_{0}$ is of type $II^{*}$ and 
$X_{\infty}$ is of type $IV^{*}$.
First determine $\text{ord}(\overline{g})$. In apriori 
$\text{ord}(\overline{g}) = 1$, $3$ or $9$. If $\text{ord}(\overline{g}) = 1$, 
a smooth fiber $E$ is $g-$stable and 
$(g\vert E)^{*}\omega_{E} = \zeta_{9}\omega_{E}$. However there exists 
no such elliptic curve. If $\text{ord}(\overline{g}) = 9$, then 
$\overline{g}$ permutes nine fibers $\{X_{\zeta_{9}^{i}t}\}_{i = 0}^{p-1}$, 
and there exists an integer $m$ with 
$24 = \chi_{top}(X) = \chi_{top}(X_{0}) + \chi_{top}(X_{\infty}) + 9m 
= 18 + 9m $, a contradiction. Thus, $\text{ord}(\overline{g}) = 3$. 
Then $g^{3}$ acts on each fiber $(g\vert E)^{*}\omega_{E} = 
\zeta_{3}\omega_{E}$. Thus, the $J-$invariant map
$J : {\bold P}^{1} \rightarrow {\bold P}^{1}$ is constant 
$j({\bold C}/ {\bold Z} + {\bold Z} \zeta_{3}) =0$. 
In particular, each singular 
fiber is either of Type $II$, $II^{*}$, $IV$ or $IV^{*}$.
%by the classification of singular fibers. 
Then by counting Euler number of $\chi_{top}(X)$, 
we see that there exists three other singular fibers of $\Phi$ of type $II$ 
permuted by $g$. Thus, we may adjust an inhomogeneous coordinate $t$  so that 
singular fibers of $\Phi$ are $X_{0}$, $X_{\infty}$ and $X_{\zeta_{3}^{i}}$ 
($i = 0, 1, 2$). Now by the same argument as before, we can readily see that 
$X$ is isomorphic to the Jacobian K3 surface 
$y^{2} = x^{3} +t^{5}(t^{3} -1)$. 
\par 
Finally consider the case where $N = 27$. 
As in the previous case, we readily see that 
$\text{ord}(\overline{g}) = 9$, the $J-$invariant map is the constant map
$J(t) = j({\bold C}/ {\bold Z} + {\bold Z} \zeta_{3}) =0$, the reducible 
singular fiber is of Type $IV$ and that the remaining singular fibers consist 
of one singular fiber of Type $II$ stable under $g$ and nine singular fibers 
of Type $II$ permuted by $g$. Then, we may normalise inhomogeneous 
coordinate $t$ of the base so that  $X_{0}$ and $X_{\zeta_{9}^{i}}$ 
($0 \leq i \leq 8$) are of Type $II$ and $X_{\infty}$ is of Type $IV$. 
Now, writing Weierstrass equation and adjusting coordinates of fibers suitably
just as before, we can readily see that $X$ is isomorphic to the Jacobian K3 
surface 
$y^{2} = x^{3} +t(t^{9} -1)$. 

\par
This complete the uniqueness for the case where $N = 3, 3^{2}$, or $3^{3}$. 

\par \vskip 2pc
{\bf \S 3. Determination of singular fibers when  $h(X)$  
equals a prime  $p$ ($\ge 5$)  and satisfies  $\varphi(h(X)) = t(X)$}

\par \vskip 0.5pc 
Let $p \geq 5$ be a prime number in $\Omega$ and 
$X$ a K3 surface with $\varphi(h(X)) = t(X)$ and with $h(X) = p$. 
Let us fix a solution of $4t^{p} + 27 =0$ and denot it by 
$\alpha_{p}$.
\par \vskip 1pc
The goal of this section is to show the following:
\par \vskip 1pc 

{\bf Proposition (3.1).}
{\it For each $p$, $X$ admits a Jacobian fibration 
$\Phi_{p} : X \rightarrow \Bbb P^{1}$ whose singular fibers are as follows:

\par \vskip 0.5pc
$X_{0}$  is of Type $II$, $X_{\infty}$  is of Type $III$, and 
$X_{\alpha_{19}\zeta_{19}^{i}}$ ($1 \leq i \leq 19$) is of Type $I_{1}$ 
when  $p = 19$;

\par 
$X_{0}$  is of Type $IV$, $X_{\infty}$  is of Type $III$, and 
$X_{\alpha_{17}\zeta_{17}^{i}}$ ($1 \leq i \leq 17$) is of Type $I_{1}$ 
in the case where $p = 17$;

\par 
$X_{0}$ is of Type $II$, $X_{\infty}$ is of Type $III^{*}$, and 
$X_{\alpha_{13}\zeta_{13}^{i}}$ ($1 \leq i \leq 13$) is of Type $I_{1}$ 
in the case where $p = 13$;
\par 

$X_{0}$ is of Type $II^{*}$, $X_{\infty}$ is of Type $III$, and 
$X_{\alpha_{11}\zeta_{11}^{i}}$ ($1 \leq i \leq 11$) is of Type $I_{1}$
in the case where $p = 11$;

\par
$X_{0}$ is of Type $IV^{*}$, $X_{\infty}$ is of Type $III^{*}$, and 
$X_{\alpha_{7}\zeta_{7}^{i}}$ ($1 \leq i \leq 7$) is of Type $I_{1}$
in the case where $p = 7$;

\par
$X_{0}$ is of Type $II^{*}$, $X_{\infty}$ is of Type $III^{*}$, and 
$X_{\alpha_{5}\zeta_{5}^{i}}$ ($1 \leq i \leq 5$) is of Type $I_{1}$ 
in the case where $p = 5$.}

\par \vskip 1pc 
{\it Proof.} By (1.1), there is a Jacobian fibration 
$\Phi : X \rightarrow {\bold P}^1$.
For a generator  $g$  of  $H_X$, there is an element  
${\overline g} \in$ Aut$({\bold P}^1)$
such that  ${\overline g} \circ \Phi = \Phi \circ g$
because  $g^*|S_X = id$. Note also that each smooth rational curve 
on $X$ is $g-$stable.

\par\vskip 1pc 

{\bf Claim(3.2).} ${\overline g}$  is of order  $p$.

\par\vskip 1pc

{\it Proof.} Suppose the contrary that the assertion is false. 
Then  ${\overline g} = id$.  Let  $E$  be a smooth
fiber of  $\Phi$.  Then  $g(E) = E$.
Since  $\omega_E \wedge \Phi^*(dt)$  gives a
nowhere vanishing 2-form around  $E$, $g^* \omega =
\zeta_p \omega$  implies that
$(g|E)^* \omega_E = \zeta_p \omega_E$.
But there is no such elliptic curve
with such action. 

\par \vskip 1pc
We adjust an inhomogeneous coordinate  $t$  of
${\bold P}^1$  such that  $({\bold P}^1)^{\overline g}
= \{0, \infty\}$.  Then only $X_0$ and $X_{\infty}$
are the $g$-stable fibers.
Note that singular fibers
$X_a$  where  $a \ne 0, \infty$ (and hence  $X_a$ 
is not  $g$-stable) is of Kodaira type
$I_1$  or  $II$, for otherwise  $X_a$  contains
a smooth rational curve which is  $g$-stable
for  $g^*|S_X = id$. Since ${\overline g}$  permutes  
$\{X_a, X_{\zeta_p a}, \dots, X_{\zeta_p^{p-1}a} \}$, we have
$$(3.0.1) \,\,\,\,\, 24 = \chi_{top}(X_0) + \chi_{top}(X_{\infty}) +
p c_1 + 2p c_2,$$  where
$p c_1, p c_2$  denote the numbers of singular fibers
of types $I_1, II$, respectively.
Moreover, $X^g = (X_0)^g \coprod (X_{\infty})^g$, whence  
$$(3.0.2) \,\,\,\,\, \chi_{top}(X^g) = \chi_{top}((X_0)^g) +
\chi_{top}((X_{\infty})^g).$$

\par \vskip 1pc
{\bf Lemma(3.3).} {\it When  $X_t$  is smooth (i.e.,
of type  $I_0$), we set  $n_t = 0$, and when  
$X_t$  is singular, we let  $n_t$  denote
the number of irreducible components of  $X_t$.
Then each of  $X_0$  and  $X_{\infty}$
is either of type  $I_{pm}, I_{pm}^*, II, III, IV, II^*, III^*,
IV^*$.  For both  $t = 0, \infty$,
$\chi_{top}(X_t) = \chi_{top}((X_t)^g) = 
n_t$ (resp. $n_t + 1$) if  $X_t$  is of type
$I_{pm}$ (resp. otherwise).}

\par \vskip 1pc
{\it Proof.} We only consider  $X_0$, for  $X_{\infty}$
is exactly the same.

\par
By the classification of elliptic fibers, $\chi_{top}(X_0) 
= n_0$ (resp. $n_0 + 1$)  if  $X_0$  is of type  $I_{n_0}$
(resp. otherwise).  We now show that  $\chi_{top}(X_0^g)
= \chi_{top}(X_0)$.  

\par
If  $X_0$  is a smooth fiber then either  $X_0 \subseteq X^g$
or  $X_0 \cap X^g = \phi$  because there is no elliptic curve
with an automorphism  $g$  of prime order  $p$ ($\ge 5$)
fixing at least one point.  It follows that
$\chi_{top}(X_0^g) = \chi_{top}(X_0) = 0 = n_0$  in 
this case.

\par
Now assume that  $X_0$  is singular. Notice the following facts 
(cf. 3-Go lemma in [OZ1, \S 2]): 
\par \vskip 0.5pc 
(1) If $Q \in X_0^g$, then there exist a local coordinates $(x_{Q}, y_{Q})$ 
around $Q$ and an integer $a$ such that $g^{*}(x_{Q}, y_{Q}) = 
(\zeta_{p}^{a}x_{Q}, \zeta_{p}^{-a+1}y_{Q})$ (as 
$g^{*}\omega_{X} = \zeta_{p}\omega_{X}$);
\par 
(2) If $g \vert C \not= id$ for a smooth rational curve $C$, then 
$C^{g}$ consists of two points, say, $Q_{1}$, $Q_{2}$. If 
$(g\vert C)^{*}(t_{Q_{1}}) = \zeta_{p}^{b}t_{Q_{1}}$ around $Q_{1}$, then 
$(g\vert C)^{*}(t_{Q_{2}}) = \zeta_{p}^{-b}t_{Q_{2}}$ around $Q_{2}$. 
\par \vskip 0.5pc 
Now, using these facts and passing to the normalisation of $X_{0}$ in 
the case of Types $I_{1}$ and $II$, we can identify  $X_0^g$  for each possible
type of  $X_0$  and hence deduce easily the result. 
 
\par \vskip 1pc
{\bf Claim(3.4).}  {\it We have  
$\chi_{top}(X_0) + \chi_{top}(X_{\infty}) = 24 - p$.
In particular, all singular fibers other than  $X_0, X_{\infty}$,
are of type  $I_1$. Moerover these are permuted by $g$.}  

\par \vskip 1pc
{\it Proof.} By (3.0.2) and (3.2),
$$\chi_{top}(X_0) + \chi_{top}(X_{\infty})
= \chi_{top}(X_0^g) + \chi_{top}(X_{\infty}^g) =
\chi(X^g) = $$
$$\sum_{i=0}^4 tr(g^*|H^i(X, {\Bbb Z}))
= 2 + tr(g^*|S_X) + tr(g^*|T_X) =$$ 
$$2 + (22-(p-1)) + (-1) = 24 - p.$$

\par
Now (3.0.1) implies that  $24 = \chi_{top}(X) = (24-p) + p c_1 + 2p c_2$,
and  $c_1 + 2c_2 = 1$.  Hence  $c_1 = 1, c_2 = 0$.
This proves Claim(3.4).

\par \vskip 1pc
{\bf Lemma(3.5).} {\it The pair of  $g$-stable fibers  $(X_0, X_{\infty})$  
of the elliptic fibration  $\Phi : X \rightarrow {\bold P}^1$  
is one of the following types, after switching the indices  
$0, \infty$  if necessary:}

\par \vskip 0.5pc
{\it $(II, III)$  {\it if}  $p = 19$;

\par
$(IV, III)$  if  $p = 17$; 

\par
$(II, III^*)$, or  $(IV^*, III)$  if  $p = 13$; 

\par
$(II^*, III)$, or  $(IV, III^*)$, or  $(I_{11}, II)$  if  $p = 11$; 

\par
$(IV^*, III^*)$, or  $(IV, I_7^*)$, or
$(I_7, II^*)$, or  $(III, I_{14})$  if  $p = 7$; 

\par
$(II^*, III^*)$, or  $(IV^*, I_5^*)$, 
$(III, I_{10}^*)$, $(III^*, I_{10})$, or  $(IV, I_{15})$  
{\it if}  $p = 5$.} 

\par \vskip 1pc
{\it Proof.} This readily follows from (3.3) and (3.4).
\par \vskip 1pc 
In order to complete (3.1), it is enough to show the following:

\par \vskip 1pc
{\bf Lemma(3.6).} {\it In Lemma $(3.5)$, replacing  $\Phi$
by a new one,  we may assume that  $(X_0, X_{\infty})$
has the following type:
$(II, III)$, or  $(IV, III)$, or  $(II, III^*)$,
or  $(II^*, III)$, or  $(IV^*, III^*)$, or
$(II^*, III^*)$  if  $p = 19$, or  $17$, or  $13$, or  $11$, or
$7$  or  $5$.}

\par \vskip 1pc
{\it Proof.}  In the case  $p = 5$ (resp. $p = 7$  or  $p = 11$),
$(X_0, X_{\infty})$  has one of 
5 (resp. 4, 3) types in (3.5).
Suppose that  $(X_0, X_{\infty})$  is not of the first 
type in (3.5).  Let  $F$  be a section of  $\Phi$.  Clearly,
$X_0 + F + X_{\infty}$  contains a weighted rational tree
$X_0''$  of Kodaira type  $II^*$  (resp. $III^*$, or  $II^*$).
Then  $X_0''$  is nef.
Now the Riemann-Roch theorem implies that
there is an elliptic fibration  $\Psi$  on
$X$  with  $X_0''$  as a ($g$-stable) fiber.

\par
It is easy to see that  $X_0$  or  $X_{\infty}$  
contains a cross-section of  $\Psi$.  
Applying  (3.4) to  $\Psi$, we see that the only two  
$g$-stable fibers of  $\Psi$  are of the first type
in (3.4).  Now (3.6) follows by replacing  $\Phi$  by  $\Psi$.

\par\vskip 1pc 
Next consider the case $p = 13$.  Suppose that the pair of the only two  $g$-stable
fibers  $(X_0, X_{\infty})$  is of the second type
$(IV^*, III)$  in (3.5).  

\par \vskip 1pc
{\bf Claim.} There are two cross-sections  $F_1, F_2$  of
$\Phi$  such that  $F_1 \cap F_2 = \phi$  and that
$F_1$  and  $F_2$  meet different (multiplicity one)
components in  $X_t$  for both  $t = 0, \infty$.

\par \vskip 1pc
Once Claim is proved to be true, (3.6) follows
by replacing  $\Phi$  by the elliptic fibration
one of whose singular fibers is of type  $III^*$
and contained in  $X_0 + F_1 + F_2 + X_{\infty}$.

\par \vskip 1pc
Now we prove Claim.
We fully use the notation and results in 
[Sh, Theorems 8.4, 8.6 and 8.7].
Fix one section  $F_1$  as the zero in the Mordell-Weil
lattice  $E(K)$  of  $\Phi$.  First, $E(K)$  is torsion free.
Indeed, if  $F$ ($\ne F_1$)  is a torsion in  $E(K)$,
then the height pairing  $0 = <F_1, F_1> =
2 \chi({\Cal O}_X) + 2F . F_1 - \sum_{v \in R} contr_v(F)
= 4 + 2F . F_1 -$ ($4/3$  or  0) $-$ ($1/2$  or  0) 
$\ge 2F . F_1 + 13/6 \ge 13/6 > 0$, a contradiction.
So  $E(K)$  is a torsion free lattice of rank 1 [Sh, Cor. 5.3].
Write  $E(K) = {\Bbb Z} F_2$. 

\par
Denoting by  $n$  the index of the sublattice  $E(K)^0$
in  $E(K)$, we have  
$n^2<F_2, F_2> =\text{det}(E(K)^0) = (\text{det} S_X) n^2/(3 \times 2)$,
and  $<F_2, F_2> = 13/6$.
Now the equality  $13/6 = <F_2, F_2> =
2 \chi({\Cal O}_X) + 2F_2 . F_1 - \sum_{v \in R} contr_v(F_2)$
and the description of  $contr_{v}(F_2)$  in [S, (8.16)]
imply Claim.
This also completes the proof of (3.6). 

\par \vskip 2pc
{\bf \S 4. Weierstrass equations of K3 surfaces when
$h(X)$  equals a prime  $p$ ($\ge 5$)  and satisfies
$\varphi(h(X)) = t(X)$}

\par \vskip 0.5pc
Let $y^{2} = x^{3} + a_{p}(t)x + b_{p}(t)$ be the minimal Weierstrass 
equation of $\Phi_{p} : X \rightarrow \Bbb P^{1}$ in (3.1). In this sectin, 
we determine this equation for each $p$ by applying the N\'eron-Tate algorithm
([Ne, Table in the last page]). This will implies the 
uniqueness of a K3 surface $X$ with $\varphi(h(X)) = t(X)$ and with 
$h(X) = p \geq 5$ for each $p$.
\par \vskip 1pc

Since $g$ acts on the base as 
$\overline{g}^{*}(t) = \zeta_{p}^{k}t$ 
(for some $k$ with $(k, p) = 1$), 
the $J-$ invariant function $J_{p}(t) := 4 a_{p}(t)^{3}/\Delta_{p}(t)$ is 
$<\zeta_p>-$invariant, and 
$\Delta_{p}(t) :=  4 a_{p}(t)^{3} + 27 b_{p}(t)^{2}$, which defines 
the discriminant divisor of $\Phi_{p}$, is semi $<\zeta_p>-$invariant. 
Thus, $a_p(t)$ is semi $<\zeta_p>-$invariant. Since $J_p(t) \not= 0$, 
we have $a_{p}(t) \not= 0$. This together with the invariance of $J_p(t)$ 
also implies the semi-invariance of $b_p(t)$.

\par
On the other hand, by the description of singular fibers and by the fact that 
$\text{deg}\Delta_{p}(t) \leq 24$, 
we have 
$\Delta_{19}(t) = C_{19}t^{2}(4t^{19}+27)$; 
$\Delta_{17}(t) = C_{17}t^{4}(4t^{17}+27)$; 
$\Delta_{13}(t) = C_{13}t^{2}(4t^{13}+27)$; 
$\Delta_{11}(t) = C_{11}t^{10}(4t^{11}+27)$;
$\Delta_{7}(t) = C_{7}t^{8}(4t^{7}+27)$; 
$\Delta_{5}(t) = C_{5}t^{10}(4t^{5}+27)$.
Here $C_{p} \not= 0$ are some constants.
Moreover, in each case, the singular fiber $X_{\infty}$ is the form of 
the finite quotient of $\Bbb C/(\Bbb Z + \Bbb Z\zeta_{4})$. 
Then we have $1 = J_{p}(\infty) = lim_{t \rightarrow \infty} J_{p}(t)$. 
This implies:
$a_{p}(t) = A_{p}t^{7}$ if $p = 19, 17, 11$, 
$a_{p}(t) = A_{p}t^{5}$ if $p = 13, 7$ and (usinig also the semi-invariance)
$a_{5}(t) = A_{5}t^{5} + C$ if $p = 5$, 
where $A_{p}$ are constants with $A_{p}^{3} = C_{p}$. 
In the case $p=5$, using
$\Delta_{5}(t) = 4a_{5}(t)^{3} + 27b_{5}(t)^{2}$ and 
the semi-invariance of $b_{5}(t)$, we readily see that $C= 0$. 
Thus, $a_{5}(t) = A_{5}t^{5}$. 

\par
Now, substituting these into $\Delta_{p}(t) = 4a_{p}(t)^{3} + 27b_{p}(t)^{2}$,
we obtain 
$b_{19}(t) = B_{19}t$;
$b_{17}(t) = B_{17}t^{2}$;
$b_{13}(t) = B_{13}t$;
$b_{11}(t) = B_{11}t^{5}$; 
$b_{7}(t) = B_{7}t^{4}$;
$b_{5}(t) = B_{5}t^{5}$, 
where $B_{p}$ are constants with $B_{p}^{2} = C_{p}$.
Then, there exists a constant $D_{p} \not= 0$ such that 
$A_{p} = D_{p}^{4}$ and $B_{p} = D_{p}^{6}$. 
Thus, the Weierstrass equation of $X_{p}$ is:

\par
$y^{2} = x^{3} + D_{p}^{4}t^{7}x + D_{p}^{6}t$ if $p = 19$;
$y^{2} = x^{3} + D_{p}^{4}t^{7}x + D_{p}^{6}t^{2}$ if $p = 17$;
$y^{2} = x^{3} + D_{p}^{4}t^{5}x + D_{p}^{6}t$ if $p = 13$;
$y^{2} = x^{3} + D_{p}^{4}t^{7}x + D_{p}^{6}t^{5}$ if $p = 11$;
$y^{2} = x^{3} + D_{p}^{4}t^{5}x + D_{p}^{6}t^{4}$ if $p = 7$;
$y^{2} = x^{3} + D_{p}^{4}t^{5}x + D_{p}^{6}t^{5}$ if $p = 5$.
\par 

Now changing the coordinates of fibers $(x, y)$ by $(D_{p}^{2}x, D_{p}^{3}y)$, 
we can normalise the equation as:

\par 
$y^{2} = x^{3} + t^{7}x + t$ if $p = 19$;
$y^{2} = x^{3} + t^{7}x + t^{2}$ if $p = 17$;
$y^{2} = x^{3} + t^{5}x + t$ if $p = 13$;
$y^{2} = x^{3} + t^{7}x + t^{5}$ if $p = 11$;
$y^{2} = x^{3} + t^{5}x + t^{4}$ if $p = 7$;
$y^{2} = x^{3} + t^{5}x + t^{5}$ if $p = 5$.

\par 
This shows the uniqueness of a K3 surface $X$ with 
$\varphi(h(X)) = t(X)$ and with $h(X) = p \geq 5$ for each $p$. 

\par \vskip 1pc
{\bf \S 5. Conclusion.}
\par \vskip 0.5pc 
In this section, we completes the proof of the uniqueness part of 
Theorem 2(2) and Corollaries 3 and 4.
\par \vskip 0.5pc 
The uniqueness part of Theorem 2(2) follows from Section 2 (the case where 
$h(X) = 3, 3^{2}, 3^{3}$), Section 4 
(the case where $h(X) =p \geq 5$ is prime) and [MO, Theorem3] 
(the case where $h(X) = 5^{2}$). Q.E.D.

\par \vskip 0.5pc 
Next we show Corollary 3. Set $p = 19$ (resp. $17$ or $13$). 
Since $g^{*}\omega_{X} \not= \omega_{X}$ by [Ni1, \S 5], $g^{*}\vert T_{X}$ 
is of order $p$.
Then $t(X) = p-1$ by [Ni1, Theorem 3.1 and Cor. 3.3] whence 
$\text{rank}S_{X} = 22 - (p-1) = 4$ (resp. $6$ or $10$). In each case, 
$\text{rank}S_{X} < \varphi(p) = p-1$. This implies $g^{*} \vert S_{X} = id$, 
whence $<g> \subset H_{X}$. Combining this with 
Theorems 1(1) and 2(1), we get $H_{X} = <g>$. Now we may apply Theorem 2(2) 
to conculde the result. Q.E.D.
\par \vskip 0.5pc 
Finally, we show Corollary 4. Let $\overline{X}$ be the canonical cover of $Z$,
$<g>$ the Galois group of this covering and $X$ the minimal resolution of 
$\overline{X}$. Then $X$ is a K3 surface and $g$ 
induces an automorphism of $X$ of order $I(Z)$. Now we can apply Corollary 3 
to get $(X, <g>) \simeq (X_{I}, <g_{I}>)$. Since $\overline{X} \rightarrow Z$ 
has no ramification curves, every $g-$fixed curve on $X$ must be contracted 
under $X \rightarrow \overline{X}$.  Now the result follows from 
the maximality assumption (*) on $Z$. Q.E.D.
  
\par 

\par \vskip 4pc
Keiji Oguiso
\par
Department of Mathematical Sciences
\par
University of Tokyo
\par
Komaba, Meguro, Tokyo
\par
JAPAN
\par
e-mail: oguiso$\@$ms.u-tokyo.ac.jp

\par \vskip 1pc
De-Qi Zhang
\par
Department of Mathematics
\par
National University of Singapore
\par
Lower Kent Ridge Road 
\par
SINGAPORE 119260
\par
e-mail: matzdq$\@$math.nus.edu.sg

\enddocument